\documentclass[12pt, leqno]{amsart}

\usepackage{texdraw,multicol,rotating}
\usepackage{color}
\usepackage{enumerate}

\def\oh{\overline{h}}

\usepackage{amsmath}
\usepackage{amssymb}
\usepackage{enumerate}
\usepackage[all]{xy}

\input{txdtools.tex}

\numberwithin{equation}{section}

\theoremstyle{plain}
\newtheorem{thm}{Theorem}[section]
\newtheorem{prop}[thm]{Proposition}
\newtheorem{lem}[thm]{Lemma}
\newtheorem{cor}[thm]{Corollary}

\theoremstyle{definition}
\newtheorem{defi}[thm]{Definition}

\newtheorem{example}[thm]{Example}

\theoremstyle{remark}

\newbox{\tmpa}
\newbox{\tmpb}

\newcommand{\nc}{\newcommand}
\nc{\Uq}{U_q} \nc{\Z}{\mathbf{Z}} \nc{\C}{\mathbf{C}}
\nc{\Q}{\mathbf{Q}}
\nc{\op}{\oplus} \nc{\ot}{\otimes} \nc{\pv}{P^{\vee}}
\nc{\ali}{\alpha_i} \nc{\B}{\mathbf{B}} \nc{\F}{\mathbf{F}}
\nc{\bP}{\mathbf{P}} \nc{\V}{\mathbf{V}} \nc{\La}{\Lambda}
\nc{\la}{\lambda} \nc{\nbinom}[2]{\genfrac{}{}{0pt}{1}{{#1}}{{#2}}}
\nc{\qbinom}[2]{\left[\genfrac{}{}{0pt}{1}{{#1}}{{#2}}\right]}
\nc{\path}{\mathcal{P}} \nc{\fit}{\tilde{f}_i}
\nc{\eit}{\tilde{e}_i} \nc{\fjt}{\tilde{f}_j} \nc{\ejt}{\tilde{e}_j}
\nc{\Y}{\mathbf{Y}} \nc{\A}{\mathbf{A}} \nc{\ra}{\rightarrow}
\nc{\vep}{\varepsilon} \nc{\vphi}{\varphi} \nc{\vp}{\varphi}
\nc{\g}{\mathfrak{g}} \nc{\h}{\mathfrak{h}} \nc{\oP}{\overline{P}}
\nc{\pathp}{\mathbf{p}} \nc{\N}{\mathcal{N}} \nc{\R}{\mathcal{R}}
\nc{\loops}{{\mathrm{loop}}}
\nc{\gl}{\mathfrak{gl}}
\nc{\Tr}{\on{Tr}}
\nc{\End}{\on{End}}
\nc{\epsor}{\varepsilon^{\mathrm{or}}}
\nc{\reg}{\mathrm{reg}}
\nc{\st}{\mathrm{st}}
\nc{\Ker}{\on{Ker}}

\nc{\tris}{ \bsegment \move(0 0)\lvec(10 0)\lvec(10 10)\lvec(0
0)\ifill f:0.7 \esegment } \nc{\recs}{ \bsegment \move(0 0)\lvec(10
0)\lvec(10 5)\lvec(0 5)\lvec(0 0)\ifill f:0.7 \esegment }
\nc{\hcvec}[5]{%
\getpos(#1 #3)\spx\spy \getpos(#2 #3)\epx\epy \getpos(#4
#5)\xoff\yoff \realadd \spx \xoff \twox \realadd \epx {-\xoff} \thrx
\realadd \spy \yoff \posy \move({\spx} {\spy}) \clvec ({\twox}
{\posy})({\thrx} {\posy})({\epx} {\epy}) \rmove(0 0) }

\nc{\ahead}[2]{%
\cossin (0 0)({#1} {#2})\cosa\sina \bsegment
  \drawdim in \setunitscale 0.065
  \realmult {-0.5} \cosa \hcosa
  \realmult {-0.5} \sina \hsina
  \move({\hcosa} {\hsina}) \ravec({\cosa} {\sina})
\esegment }
\nc{\boxi}{%
{%
\savebox{\tmppic}{\begin{texdraw} \small \drawdim em \textref h:C
v:C \setunitscale 0.55 \htext(0 0){$i$} \move(-1 -1)\lvec(-1
1)\lvec(1 1)\lvec(1 -1)\lvec(-1 -1)
\end{texdraw}}%
\raisebox{-0.19\height}{\usebox{\tmppic}}%
}%
}
\nc{\boxj}{%
{%
\savebox{\tmppic}{\begin{texdraw} \small \drawdim em \textref h:C
v:C \setunitscale 0.55 \htext(0 0.1){$j$} \move(-1 -1)\lvec(-1
1)\lvec(1 1)\lvec(1 -1)\lvec(-1 -1)
\end{texdraw}}%
\raisebox{-0.19\height}{\usebox{\tmppic}}%
}%
}
\nc{\boxipo}{%
{%
\savebox{\tmppic}{\begin{texdraw} \small \drawdim em \textref h:C
v:C \setunitscale 0.55 \htext(0.15 0){$i\!\!+\!\!1$} \move(-1.4
-1)\lvec(-1.4 1)\lvec(1.4 1)\lvec(1.4 -1)\lvec(-1.4 -1)
\end{texdraw}}%
\raisebox{-0.19\height}{\usebox{\tmppic}}%
}%
} \everytexdraw{ \drawdim in \arrowheadsize l:0.065 w:0.03
\arrowheadtype t:F \textref h:C v:C }
\newsavebox{\tmppic}
\newsavebox{\tmpfig}
\newsavebox{\tmpdraw}
\newsavebox{\tmpfiga}
\newsavebox{\tmpfigb}
\newsavebox{\tmpfigc}
\newsavebox{\tmpfigd}
\newsavebox{\tmpfige}
\newsavebox{\tmpfigf}
\newsavebox{\tmpfigg}
\newsavebox{\tmpfigh}
\newsavebox{\tmpfigi}
\newsavebox{\tmpfigj}
\newsavebox{\tmpfigk}
\newsavebox{\tmpfigl}
\newsavebox{\tmpfigm}
\newsavebox{\tmpfign}
\newsavebox{\tmpfigo}
\newsavebox{\tmpfigp}
\newsavebox{\tmpfigq}
\newsavebox{\tmpfigr}
\newsavebox{\tmpfigs}
\newsavebox{\tmpfigt}
\newsavebox{\tmpfigu}
\newsavebox{\tmpfigv}
\newsavebox{\tmpfigw}
\newsavebox{\tmpfigx}
\newsavebox{\tmpfigy}
\newsavebox{\tmpfigz}

\newsavebox{\tmpfigaa}
\newsavebox{\tmpfigab}
\newsavebox{\tmpfigac}
\newsavebox{\tmpfigad}
\newsavebox{\tmpfigae}
\newsavebox{\tmpfigaf}
\newsavebox{\tmpfigag}
\newsavebox{\tmpfigah}
\newsavebox{\tmpfigai}
\newsavebox{\tmpfigaj}
\newsavebox{\tmpfigak}
\newsavebox{\tmpfigal}
\newsavebox{\tmpfigam}
\newsavebox{\tmpfigan}
\newsavebox{\tmpfigao}
\newsavebox{\tmpfigap}
\newsavebox{\tmpfigaq}
\newsavebox{\tmpfigar}
\newsavebox{\tmpfigas}
\newsavebox{\tmpfigat}
\newsavebox{\tmpfigau}
\newsavebox{\tmpfigav}
\newsavebox{\tmpfigaw}
\newsavebox{\tmpfigax}
\newsavebox{\tmpfigay}
\newsavebox{\tmpfigaz}

\newsavebox{\tmpfigba}
\newsavebox{\tmpfigbb}
\newsavebox{\tmpfigbc}
\newsavebox{\tmpfigbd}
\newsavebox{\tmpfigbe}
\newsavebox{\tmpfigbf}
\newsavebox{\tmpfigbg}
\newsavebox{\tmpfigbh}

\nc{\node}{\lcir r:1 }
\nc{\sline}{\bsegment\savepos(10 0)(*ex *ey)
            \move(1 0)\rlvec(8 0)
            \esegment\move(*ex *ey)}
\nc{\dline}{\bsegment\savepos(10 0)(*ex *ey)
            \move(0.93 0.4)\rlvec(8.14 0)\rmove(0 -0.8)\rlvec(-8.14 0)
            \esegment\move(*ex *ey)}
\nc{\uline}{\bsegment\savepos(0 10)(*ex *ey)
            \move(0 1)\rlvec(0 8)
            \esegment\move(*ex *ey)}
\nc{\lpoint}{\savecurrpos(*ex *ey)
             \rmove(2.5 1.5)\rlvec(-1.5 -1.5)\rlvec(1.5 -1.5)
             \move(*ex *ey)}
\nc{\rpoint}{\savecurrpos(*ex *ey)
             \rmove(-2.5 -1.5)\rlvec(1.5 1.5)\rlvec(-1.5 1.5)
             \move(*ex *ey)}
\nc{\bline}{\bsegment\savepos(10 0)(*ex *ey)
            \linewd 0.6 \move(1.1 0)\rlvec(7.8 0)
            \esegment\move(*ex *ey)}

\nc{\araise}[1]{\raisebox{4.5pt}{#1}}
\nc{\braise}[1]{\raisebox{12.1pt}{#1}}
\nc{\craise}[1]{\raisebox{8pt}{#1}}
\nc{\draise}[1]{\raisebox{12pt}{#1}} \nc{\be}{\begin{enumerate}}
\nc{\ee}{\end{enumerate}} \nc{\bnum}{\begin{enumerate}[{\rm(i)}]}
\nc{\cl}{\colon} \nc{\seteq}{\mathbin{:=}} \nc{\re}{\mathrm{re}}
\nc{\im}{\mathrm{im}} \nc{\ran}{\rangle} \nc{\lan}{\langle}
\nc{\on}{\operatorname} \nc{\gr}{\mathrm{gr}}
\newcommand{\set}[2]{\left\{#1\,;\,#2\,\right\}}
\nc{\Hom}{\on{Hom}} \nc{\Homi}{\on{Hom^{\gr}}}
\nc{\Oint}{\mathcal{O}_{\mathrm{int}}}
\nc{\wt}{\on{wt}}
\nc{\Wt}{\on{Wt}} \nc{\pP}{\widetilde{P}} \nc{\eq}{\begin{eqnarray}}
\nc{\eneq}{\end{eqnarray}} \nc{\eqn}{\begin{eqnarray*}}
\nc{\eneqn}{\end{eqnarray*}} \nc{\Lemma}{\begin{lem}}
\nc{\enlemma}{\end{lem}}

\nc{\hs}{\hspace*} \nc{\bfi}{\mathbf{i}} \nc{\eps}{\varepsilon}
\nc{\ba}{\begin{array}} \nc{\ea}{\end{array}}
\renewcommand{\Im}{\operatorname{Im}}
\nc{\tf}{\tilde{f}} \nc{\id}{\operatorname{id}} \nc{\bl}{\bigl}
\nc{\br}{\bigr} 
\nc{\Irr}{\on{Irr}}
\nc{\Id}{\on{Id}}
\nc{\out}{\operatorname{out}}
\nc{\sink}{\operatorname{in}}
\nc{\ol}{\overline}

\begin{document}

\title[Geometric Construction of Highest Weight Crystals]
      {Geometric Construction of Highest Weight Crystals \\for Quantum Generalized Kac-Moody Algebras}
\author[S.-J. Kang, M. Kashiwara, O. Schiffmann]{Seok-Jin Kang$^{1}$,
Masaki Kashiwara$^{2}$, Olivier Schiffmann}

\address{Department of Mathematical Sciences
         and
         Research Institute of Mathematics \\
         Seoul National University \\ San 56-1 Sillim-dong, Gwanak-gu \\ Seoul 151-747, Korea}

         \email{sjkang@math.snu.ac.kr}

\address{Research Institute for Mathematical Sciences \\
         Kyoto University \\ Kitashirakawa, Sakyo-Ku \\ Kyoto 606-8502, Japan}

         \email{masaki@kurims.kyoto-u.ac.jp}

\address{Universit\'e Pierre et Marie Curie \\
         D\'epartement de Math\'ematiques\\ 175 rue du Chevaleret  \\ 75013 Paris,  France}

         \email{olive@math.jussieu.fr}

\thanks{$^{1}$This research was supported by KRF Grant \# 2007-341-C00001.}
\thanks{$^{2}$This research was partially supported by Grant-in-Aid for Scientific Research (B)
18340007, Japan Society for the Promotion of Science.}

\begin{abstract}

We present a geometric construction of highest weight crystals
$B(\lambda)$ for quantum generalized Kac-Moody algebras. It is given
in terms of the irreducible components of certain Lagrangian subvarieties of
Nakajima's quiver varieties associated to quivers with edge loops.

\end{abstract}

\maketitle

\vskip 1cm


\section*{Introduction}

The 1990's saw a great deal of interesting interplay
between the geometry of quiver varieties and the representation
theory of quantum groups. One of the most exciting developments in
this direction may be Lusztig's geometric construction of {\it
canonical bases}. For a Kac-Moody algebra $\g$, he constructed a
natural basis $\B$ of the negative part of the quantum group
$U_q(\g)$ in terms of simple perverse sheaves on quiver varieties
\cite{Lus90}. The basis $\B$ yields all other canonical bases of
integrable highest weight modules through natural projections.

Around the same time, Kashiwara took an algebraic approach to
construct {\it global bases} and showed how to obtain, by passing to the
crystal limit $q=0$, {\it crystals bases} which contain
most of the combinatorial information on $U_q(\g)$ and their
integrable highest weight representations \cite{Kas91}. We denote by
$B(\infty)$ and $B(\la)$ the crystal bases of $U_q^{-}(\g)$ and
$V(\la)$, respectively. It later turned out that canonical bases
and global bases coincide \cite{GL}.

In \cite{KS97}, Kashiwara and Saito gave a geometric construction of
$B(\infty)$: the crystal $B(\infty)$ can be identified with the set
of irreducible components of Lusztig's nilpotent quiver varieties which are
certain Lagrangian subvarieties of the cotangent space to the representation
varieties of a quiver. This work was generalized by Saito to a geometric construction of
$B(\la)$ using Nakajima's quiver varieties \cite{Sai}.

For generalized Kac-Moody algebras, which were introduced by
Borcherds in his study of Monstrous Moonshine \cite{Bo88}, the
crystal basis theory was developed in \cite{JKK} and it was proved
that there exist unique crystal bases $B(\infty)$ and $B(\la)$ for
$U_q^{-}(\g)$ and $V(\la)$, respectively. In \cite{JKKS}, the notion
of {\it abstract crystals} was put forward and the authors gave some
combinatorial characterizations of $B(\infty)$ and $B(\la)$. In
\cite{KKS08}, we gave a geometric construction of $B(\infty)$ for
quantum generalized Kac-Moody algebras in terms of irreducible
components of Lusztig's quiver varieties, associated this time to
quivers which may have loop edges. The main difficulty of
this work lies in that typical simple objects sitting at a vertex
with loops may have non-vanishing self extensions. This difficulty was overcome
by requiring that certain arrows are regular semisimple.

In this article, we continue to investigate the deep connection
between the geometry of quiver varieties and the representation
theory of quantum groups. In particular, we present a geometric
construction of highest weight crystals $B(\la)$ for quantum
generalized Kac-Moody algebras . We first define certain Lagrangian
subvarieties of Nakajima's quiver varieties by imposing stability
conditions on Lusztig's quiver varieties, and consider the set
${\mathcal B}^{\la}$ of irreducible components of these Lagrangian
subvarieties. We then define the Kashiwara operators on ${\mathcal
B}^{\la}$ using generic fibrations between irreducible components so
that ${\mathcal B}^{\la}$ becomes an abstract crystal. Finally, we
show that ${\mathcal B}^{\la}$ satisfies all the properties
characterizing $B(\la)$, from which we conclude that the crystal
${\mathcal B}^{\la}$ is isomorphic to $B(\la)$.

\vskip 1cm


\section{The Crystal $B(\la)$}

\vskip 3mm

In this section, we recall the definition and basic properties of
quantum generalized Kac-Moody algebras, integrable highest weight
modules and their crystals. Let $I$ be a finite or countably
infinite index set. A {\it symmetric even integral Borcherds-Cartan
matrix} is a square matrix $A=(a_{ij})_{i, j\in I}$ such that \ (i)
$a_{ii} \in \{2, 0, -2, -4, \ldots\}$ for all $i\in I$, \ (ii)
$a_{ij} = a_{ji} \in \Z_{\le 0}$ for $i \neq j$. Let
$I^{\re}=\set{i\in I}{a_{ii}=2}$ and $I^{\im}=\set{ i \in I}{a_{ii}
\le 0 }$ and call them the set of {\it real} indices and the set of
{\it imaginary} indices, respectively.

A {\it Borcherds-Cartan datum} $(A, P, \Pi, \Pi^{\vee})$ consists of

\begin{itemize}
\item[(i)] a Borcherds-Cartan matrix $A=(a_{ij})_{i, j \in I}$,

\item[(ii)] a free abelian group $P$, the {\em weight lattice},

\item[(iii)] $\Pi=\set{\alpha_i\in P}{i\in I}$, the set of {\em
simple roots},

\item[(iv)] $\Pi^{\vee} = \set{ h_i}{ i\in I }\subset
P^\vee\seteq\Hom(P,\Z)$, the set of {\em simple coroots}
\end{itemize}
satisfying the following properties:
\begin{itemize}
\item[(a)] $\langle h_i,\alpha_j \rangle = a_{ij}$ for all $i, j \in
I$,

\item[(b)] $\Pi$ is linearly independent,

\item[(c)] for any $i\in I$, there exists $\Lambda_i\in P$ such that
$\langle h_j,\Lambda_i \rangle=\delta_{ij}$ for all $j\in I$.
\end{itemize}

\noindent We denote by $P^{+}=\set{\lambda \in P}{\langle h_i,
\lambda \rangle \ge 0 \ \ \text{for all} \ i \in I}$ the set of {\it
dominant integral weights}. We also use the notation
$Q=\bigoplus_{i\in I} \Z \alpha_i$ and $Q_{+}=\sum_{i\in I} \Z_{\ge
0} \alpha_i$.

Let $q$ be an indeterminate. For $m, n \in \Z$, define
\begin{equation*}
[n] = \dfrac{q^n - q^{-n}} {q - q^{-1}}, \qquad [n]! =
\prod_{k=1}^n [k], \qquad \left[\begin{matrix} m \\
n \end{matrix}\right] = \dfrac{[m]!}{[n]! [m-n]!}.
\end{equation*}
The {\it quantum generalized Kac-Moody algebra} $U_q(\g)$ associated
with a Borcherds-Cartan datum $(A, P, \Pi, \Pi^{\vee})$ is defined
to be the associated algebra over $\Q(q)$ with 1 generated by the
elements $e_i$, $f_i$ $(i\in I)$, $q^h$ $(h\in P^{\vee})$ subject to
the defining relations: {\allowdisplaybreaks
\begin{equation} 
\begin{aligned}
& q^0 =1, \quad q^h q^{h'} = q^{h+h'}\quad\text{for} \ \ h, h' \in
P^{\vee}, \\
& q^h e_i q^{-h} = q^{\alpha_i(h)} e_i, \quad q^h f_i q^{-h} =
q^{-\alpha_i(h)} f_i \quad\text{for $h\in P^{\vee}$, $i\in I$,} \\
& e_i f_j - f_j e_i = \delta_{ij} \dfrac{K_i - K_i^{-1}} {q -
q^{-1}} \quad\text{for $i, j \in I$, where $ K_i =
q^{h_i}$,} \\
& \sum_{k=0}^{1-a_{ij}} (-1)^k \left[\begin{matrix} 1-a_{ij} \\ k
\end{matrix} \right] e_i^{1-a_{ij}-k} e_j e_i^k =0 \quad\text{if
$i\in I^{\re}$ and $i\neq j$,}\\
& \sum_{k=0}^{1-a_{ij}} (-1)^k \left[\begin{matrix} 1-a_{ij} \\ k
\end{matrix} \right] f_i^{1-a_{ij}-k} f_j f_i^k =0 \quad
\text{if $i\in I^{\re}$ and $i\neq j$,} \\
& e_i e_j - e_j e_i = f_i f_j - f_j f_i =0 \quad\text{if
$a_{ij}=0$.}
\end{aligned}
\end{equation}
} \noindent We denote by $U_q^+(\g)$ (resp.\ $U_q^-(\g)$) the
subalgebra of $U_q(\g)$ generated by the $e_i$'s (resp.\ the
$f_i$'s).

\vskip 3mm

The following notion of {\it abstract crystals} for quantum
generalized Kac-Moody algebras was introduced in \cite{JKKS}.

\begin{defi} \label{D:abstract crystal} An {\em abstract
$U_q(\g)$-crystal} or simply a {\em crystal\/} is a set $B$ together
with the maps $\wt \cl B \rightarrow P$, $\eit, \fit \cl B \rightarrow B
\sqcup \{0\}$ and $\vep_i, \vphi_i \cl  B \rightarrow \Z \sqcup
\{-\infty\}$ $(i\in I)$ satisfying the following conditions:

\begin{itemize}
\item[(i)] $\wt(\eit b) = \wt b + \alpha_i$ if $i\in I$ and
$\eit b \neq 0$,

\item[(ii)] $\wt(\fit b) = \wt b - \alpha_i$ if $i\in I$ and
$\fit b \neq 0$,

\item[(iii)] for any $i \in I$ and $b\in B$, $\vphi_i(b) = \vep_i(b) +
\langle h_i, \wt b \rangle$,

\item[(iv)] for any $i\in I$ and $b,b'\in B$,
$\fit b = b'$ if and only if $b = \eit b'$,\label{cond7}

\item[(v)] for any $i \in I$ and $b\in B$
such that $\eit b \neq 0$, we have \be[{\rm(a)}]
\item
$\vep_i(\eit b) = \vep_i(b) - 1$, $\vphi_i(\eit b) = \vphi_i(b) + 1$
if $i\in I^\re$,
\item
$\vep_i(\eit b) = \vep_i(b)$, \ $\vphi_i(\eit b) = \vphi_i(b) +
a_{ii}$ if $i\in I^\im$, \ee

\item[(vi)] for any $i \in I$ and $b\in B$ such that $\fit b \neq 0$,
we have \be[{\rm(a)}]
\item
$\vep_i(\fit b) = \vep_i(b) + 1$, $\vphi_i(\fit b) = \vphi_i(b) - 1$
if $i\in I^\re$,
\item
$\vep_i(\fit b) = \vep_i(b)$, \ $\vphi_i(\fit b) = \vphi_i(b) -
a_{ii}$ if $i\in I^\im$, \ee

\item[(vii)] for any $i \in I$ and $b\in B$ such that $\vphi_i(b) = -\infty$, we
have $\eit b = \fit b = 0$.

\end{itemize}
\end{defi}
\noindent We will often use the notation $\wt_{i}(b) = \langle h_i,
\wt(b) \rangle$ $(i \in I, b\in B)$.

\begin{defi}\label{D:mor}
Let $B_1$ and $B_2$ be crystals.

\be[{\quad \rm(a)}]

\item A map $\psi\cl B_1\rightarrow B_2$ is a {\it crystal morphism}
if it satisfies the following properties:

\bnum
\item for $b\in B_1$, we have
\begin{equation*}
\text{$\wt(\psi(b))=\wt(b)$, $\vep_i(\psi(b))=\vep_i(b)$,
$\vp_i(\psi(b))=\vp_i(b)$ for all $i\in I$,}
\end{equation*}
\item for $b\in B_1$ and $i\in I$ with $\fit b\in B_1$, we have
$\psi(\fit b)=\fit\psi(b)$. \label{cond:crysmor2} \ee

\item A crystal morphism $\psi\cl B_1\rightarrow B_2$ is called {\it
strict} if
\begin{equation*}
\psi(\eit b)=\eit\psi(b),\,\, \psi( \fit
b)=\fit\psi(b)\quad\text{for all $i\in I$ and $b\in B_1$.}
\end{equation*}
Here, we understand $\psi(0)=0$.

\item $\psi$ is called an {\em embedding} if the underlying map
$\psi\cl B_1\rightarrow B_2$ is injective. \ee
\end{defi}

\vskip 3mm

For a pair of crystals $B_1$ and $B_2$, their {\it tensor product}
is defined to be the set
$$B_1 \ot B_2 =  \set{b_1\otimes
b_2}{b_1\in B_1, b_2\in B_2},$$ where the crystal structure is
defined as follows: \ The maps $\wt, \vep_i,\vp_i$ are given by \eqn
\wt(b\otimes b')&=&\wt(b)+\wt(b'),\\
\vep_i(b\otimes b')&=&\max(\vep_i(b),
\vep_i(b')-\wt_i(b)),\\
\vp_i(b\otimes b')&=&\max(\vp_i(b)+\wt_i(b'),\vp_i(b')). \eneqn For
$i\in I$, we define \eqn \fit(b\otimes b')&=&
\begin{cases}\fit b\otimes b'
&\text{if $\vp_i(b)>\vep_i(b')$,}\\
b\otimes \fit b' &\text{if $\vp_i(b)\le \vep_i(b')$,}
\end{cases}
\eneqn For $i\in I^\re$, we define \eqn \eit(b\otimes b')&=&
\begin{cases}\eit b\otimes b'\ &\text{if
$\vp_i(b)\ge \vep_i(b')$,}\\
b\otimes \eit b' &\text{if $\vp_i(b)< \vep_i(b')$,}
\end{cases}
\eneqn and, for $i\in I^\im$, we define \eqn \eit(b\otimes b')&=&
\begin{cases}\eit b\otimes b'\
&\text{if $\vp_i(b)>\vep_i(b')-a_{ii}$,}\\
0&\text{if $\vep_i(b')<\vp_i(b)\le\vep_i(b') -a_{ii}$,}\\
b\otimes \eit b' &\text{if $\vp_i(b)\le\vep_i(b')$.}
\end{cases}
\eneqn

\vskip 3mm

\begin{example}
Let $V(\la)$ be the irreducible highest weight $U_q(\g)$-module with
highest weight $\la \in P^{+}$. For any $i\in I$, every $v \in
V(\la)$ has a unique {\it $i$-string decomposition}
\begin{equation*}
v=\sum_{k\ge 0} f_i^{(k)} u_k, \quad \text{where} \ \ e_i u_k =0 \ \
\text{for all} \ \ k\ge 0,
\end{equation*}
and
$$f_i^{(k)}\seteq\begin{cases}
f_i^{k}/[k]!&\text{if $i$ is real,}\\
f_i^{k}&\text{if $i$ is imaginary.}
\end{cases}$$
The {\it Kashiwara operators} $\eit$, $\fit$ $(i\in I)$ are defined
by
\begin{equation*}
\eit v = \sum_{k\ge 1} f_i^{(k-1)} u_k, \qquad \fit v = \sum_{k\ge
0} f_i^{(k+1)} u_k.
\end{equation*}

Let $\A_0 = \set{f/g \in \Q(q) }{ f, g \in \Q[q], g(0) \neq 0 }$ and
let $L(\la)$ be the free $\A_0$-submodule of $V(\la)$ generated by
$$\set{\tilde f_{i_1} \cdots \tilde f_{i_r} v_{\la}} { r \ge 0,
i_k \in I },$$ where $v_{\la}$ is the highest wight vector of
$V(\la)$. Then the set
\begin{equation*}
B(\la) = \set{\tilde f_{i_1} \cdots \tilde f_{i_r} v_{\la} + q
L(\la) }{ r\ge 0, i_k \in I } \setminus \{0\} \subset L(\la)/qL(\la)
\end{equation*}
becomes a $U_q(\g)$-crystal with the maps $\wt$, $\eit, \fit$,
$\vep_i, \vp_i$ ($i\in I$) defined by
\begin{equation*}
\begin{aligned}
\wt(b) & =\la - (\alpha_{i_1} + \cdots + \alpha_{i_r}) \quad
\text{for} \ \ b=\tilde f_{i_1} \cdots \tilde f_{i_r}
v_{\la} + q L(\la), \\
\vep_i (b) & =
\begin{cases} \max \set{k\ge 0 }{ \eit^k b \neq 0 }&\text{for $i\in I^\re$,}
\\
0&\text{for $i\in I^\im$,}\end{cases}\\
\vphi_i (b) & = \vep_i(b) + \wt_i(b) \quad (i\in I).
\end{aligned}
\end{equation*}
\end{example}

\begin{example}
For each $i\in I$, we define the endomorphisms $e_i',
e_i'':U_q^{-}(\g) \rightarrow U_q^{-}(\g)$ by
\begin{equation*}
e_i u - u e_i = \dfrac{K_i e_i''(u) - K_i^{-1} e_i'(u)}{q_i -
q_i^{-1}} \qquad \text{for} \ \ u \in U_q^{-}(\g).
\end{equation*}
Then every $u\in U_q^{-}(\g)$ has a unique {\it $i$-string
decomposition}
\begin{equation*}
u=\sum_{k\ge 0} f_i^{(k)} u_k, \quad \text{where} \ \ e_i' u_k =0 \
\ \text{for all} \ \ k\ge 0.
\end{equation*}
The {\it Kashiwara operators} $\eit$, $\fit$ $(i\in I)$ are defined
by
\begin{equation*}
\eit u = \sum_{k\ge 1} f_i^{(k-1)} u_k, \qquad \fit u = \sum_{k\ge
0} f_i^{(k+1)} u_k.
\end{equation*}

Let $L(\infty)$ be the free $\A_0$-submodule of $U_q^{-}(\g)$
generated by
$$\set{\tilde f_{i_1} \cdots \tilde f_{i_r} \mathbf{1} }{ r \ge 0,
i_k \in I },$$ where $\mathbf{1}$ is the multiplicative identity in
$U_q^{-}(\g)$. Then the set
\begin{equation*}
B(\infty) = \set{\tilde f_{i_1} \cdots \tilde f_{i_r} \mathbf{1} + q
L(\infty) }{ r\ge 0, i_k \in I } \subset L(\infty)/qL(\infty)
\end{equation*}
becomes a $U_q(\g)$-crystal with the maps $\wt$, $\eit, \fit$,
$\vep_i, \vp_i$ ($i\in I$) defined by
\begin{equation*}
\begin{aligned}
\wt(b) & =- (\alpha_{i_1} + \cdots + \alpha_{i_r}) \quad \text{for}
\ \ b=\tilde f_{i_1} \cdots \tilde f_{i_r}
\mathbf{1} + q L(\infty), \\
\vep_i (b) & =
\begin{cases} \max \set{k\ge 0 }{ \eit^k b \neq 0 }&\text{for $i\in I^\re$,}
\\
0&\text{for $i\in I^\im$,}\end{cases}\\
\vphi_i (b) & = \vep_i(b) + \wt_i(b) \quad (i\in I).
\end{aligned}
\end{equation*}
\end{example}

\begin{example} \label{exam:Tla}
For $\la \in P$, let $T_{\la}=\{t_{\la} \}$ and define
\begin{equation*}
\begin{aligned}
& \wt(t_{\la})=\la, \quad \eit t_{\la} = \fit t_{\la} =0 \quad
\text{for all} \ \ i\in I, \\
& \vep_i(t_{\la}) = \vphi_i(t_{\la}) = -\infty \quad \text{for all}
\ \ i\in I.
\end{aligned}
\end{equation*}
Then $T_{\la}$ is a $U_q(\g)$-crystal.
\end{example}

\begin{example}\label{exam:C}
Let $C=\{c\}$ be the crystal with $\wt(c)=0$ and
$\vep_i(c)=\vp_i(c)=0$, $\fit c=\eit c=0$ for any $i\in I$. Then $C$
is a $U_q(\g)$-crystal isomorphic to $B(0)$. For a crystal $B$,
$b\in B$ and $i\in I$, we have
\eqn
\wt(b\otimes c)&=&\wt(b),\\
\vep_i(b\otimes c)&=&\max(\vep_i(b),-\wt_i b),\\
\vp_i(b\otimes c)&=&\max(\vp_i(b),0),\\
\eit(b\otimes c)&=&
\begin{cases}
\eit b\otimes c&\text{if $\vp_i(b)\ge0$ and $i\in I^\re$,}\\
\eit b\otimes c&\text{if $\vp_i(b)+a_{ii}>0$ and $i\in I^\im$,}\\
0&\text{otherwise,}
\end{cases}\\
\fit(b\otimes c)&=&
\begin{cases}
\fit b\otimes c&\text{if $\vp_i(b)>0$,}\\
0&\text{otherwise.}
\end{cases}
\eneqn
\noindent In general, $B \ot C$ is {\it not} isomorphic to
$B$.
\end{example}

The crystal $B(\la)$ can be characterized as follows.

\begin{prop}[\cite{JKKS}] \label{P:B(la)}
Let $\la\in P^+$ be a dominant integral weight. Then $B(\lambda)$ is
isomorphic to the connected component of $B(\infty)\otimes
T_\la\otimes C$ containing $\mathbf{1}\otimes t_\lambda\otimes c$.
\end{prop}


\section{Lusztig's Quiver variety}

\vskip 3mm

Let $(I,H)$ be a quiver. For an arrow $h \cl i\to j$ in $H$, we
write $\out(h)=i$, $\sink(h)=j$ and assume that we have
an involution $-$ of $H$
such that $\out(\ol{h})=\sink(h)$ for any $h\in H$ and that 
$-$ has no fixed point. An {\it orientation} of $H$ is a subset
$\Omega$ of $H$ such that $H=\Omega \sqcup \overline{\Omega}$. We
say that $h$ is a {\it loop} if $\out(h)=\sink(h)$. We denote by
$H^{\loops}$ the set of all loops and set $\Omega^{\loops}=\Omega \cap
H^{\loops}$.

Let $c_{ij}$ denote the number of arrows in $H$ from $i$ to $j$, and
define
\begin{equation*}
a_{ij}=\begin{cases} 2- c_{ii} = 2 - (\text{the number of loops at
$i$ in $H$})
\ \ & \text{if} \ \ i=j, \\
 - c_{ij} = -(\text{the number of arrows in $H$ from $i$ to $j$}) \
\ & \text{if} \ \ i \neq j.
\end{cases}
\end{equation*}
Then $A=(a_{ij})_{i,j \in I}$ becomes a symmetric even integral
Borcherds-Cartan matrix.

For $\alpha \in Q_{+}$, let $V=V(\alpha)=\bigoplus_{i\in I} V_{i}$
be an $I$-graded vector space with
$$\underline{\dim} \, V(\alpha) \seteq \sum_{i\in I}
(\dim V_i) \alpha_i = \alpha,$$ and let
$$X(\alpha) = \bigoplus_{h\in H}\Hom(V_{\out(h)},V_{\sink(h)}).$$
The
group $GL(\alpha):=\prod_{i\in I} GL(V_i)$ acts on $X(\alpha)$ via
$$ g \cdot x =(g_{\sink(h)} x_{h} g_{\out(h)}^{-1})_{h\in H} \quad
\text{for} \ \ g=(g_i) \in GL(\alpha), \ x=(x_{h}) \in X(\alpha).$$
The symplectic form $\omega$ on $X(\alpha)$
and the {\it moment map} $\mu=\bigl(\mu_i
\cl X(\alpha) \rightarrow \gl(V_i)\bigr)_{i\in I}$ are
given by
$$\omega(x,y)=\sum_h \epsilon(h) \Tr(x_{\oh}y_h),$$
$$\mu_i(x) = \sum_{\substack{h\in H \\ \out(h)=i}} \epsilon(h)
x_{\overline{h}} x_{h},$$ where
$$\epsilon(h)=\begin{cases} 1 & \text{if} \; h\in \Omega, \\
-1 & \text{if\;}  h\in \overline{\Omega}.
\end{cases}$$

We define {\it Lusztig's quiver variety} $\N(\alpha)$ to be the
variety consisting of all $x=(x_h)_{h\in H} \in X(\alpha)$
satisfying the following conditions:

\begin{itemize}

\item[(i)] $\mu_i(x)=0$ for all $i\in I$,

\item[(ii)] there exists an $I$-graded complete flag ${F}=(F_0 \subset
F_1 \subset F_2 \subset \cdots)$ such that
$$x_{h}(F_k) \subset
F_{k} \ \ \text{for all} \ h\in \overline{\Omega}^{\loops}, \ \
x_{h}(F_{k}) \subset F_{k-1} \ \ \text{for all} \ h\in H \setminus
\overline{\Omega}^{\loops},$$

\item[(iii)] $x_{h}$ is regular semisimple for all $h \in\ol{\Omega}^{\loops}$.

\end{itemize}

\noindent  We denote by $\Irr \N(\alpha)$ the set of irreducible
components of $\N(\alpha)$.

\vskip 3mm

Fix $i\in I$ and let $t$ be the number of loops at $i$ in $\Omega$.
Write $\Omega_i^{\loops}=\{\sigma_1, \ldots, \sigma_t\}$. Let
$\R=\C\langle x_1, \ldots, x_t, y_1, \ldots, y_t \rangle$ be the
free unital associative algebra generated by $x_i, y_i$ $(i=1, \ldots, t)$.
For $x=(x_{h})_{h\in H} \in \N(\alpha)$ and $f\in \R$, we define
\begin{equation*}
\begin{aligned}
f(x)& = f(x_{\sigma_1}, \ldots, x_{\sigma_t},
x_{\overline{\sigma_1}}, \ldots, x_{\overline{\sigma_t}})\in\End(V_i), \\
\C\langle x \rangle_{i}& = \set{ f(x)}{f \in \R }\subset\End(V_i), \\
\epsor_i (x) &= \mathrm{codim}_{V_i} \Bigl(\C \langle
x\rangle_i \cdot \sum_{\substack{h\cl j\rightarrow i \\ j \neq i}}
\Im x_{h}\Bigr).
\end{aligned}
\end{equation*}
Since $\epsor_i$ is a semicontinuous function, it takes a constant value
$\epsor_i(\Lambda)$ on an open dense subset
of any irreducible component $\Lambda$ of
$\mathcal{N}(\alpha)$.
Note that we shall see later that 
$\epsor_i(\Lambda)=\max\set{n\ge0}{\eit^n(\Lambda)\not=0}$.
We set
\eq
&&\N(\alpha)_{i,l}=\set{x \in \N(\alpha)}{\text{$ \epsor_i=l$ on a neighborhood of $x$ in $\N(\alpha)$}}.
\eneq
Then $\cup_{l\ge0}\N(\alpha)_{i,l}$ is an open dense subset of
$\N(\alpha)$.
It is shown in \cite{KKS08} that if $\Lambda \in \Irr \N(\alpha)$
and $\epsor_i(\Lambda)=0$ for all $i\in I$, then $\alpha =0$
and $\Lambda =\{0\}$.

\vskip 3mm

For each $\alpha \in Q_{+}$ and $l \ge 0$, let
\begin{equation*}
\begin{aligned}
E(\alpha, l\alpha_i)  = &\big\{(x,x', x'', \phi', \phi'') \;; \;
 x \in X(\alpha + l \alpha_i), x' \in \N(\alpha), x'' \in \N(l \alpha_i), \\
& 0 \longrightarrow V(\alpha) \overset{\phi'}\longrightarrow
V(\alpha+l\alpha_i) \overset{\phi''} \longrightarrow V(l \alpha_i)
\longrightarrow
0 \ \ \text{is exact}, \\
& \phi' \circ x' = x \circ \phi', \; \phi'' \circ x=x'' \circ
\phi''\big\}
\end{aligned}
\end{equation*}
and consider the canonical projections
\begin{equation}\label{E:P1}
X(\alpha) \times X(l\alpha_i) \overset{p_1} \longleftarrow E(\alpha,
l\alpha_i) \overset{p_2} \longrightarrow X(\alpha+l\alpha_i)
\end{equation}
 given by
\begin{equation*}
(x', x'') \leftarrow (x, x', x'', \phi', \phi'') \mapsto x.
\end{equation*}
Let $\N(\alpha, l \alpha_i) = p_2^{-1}(\N(\alpha + l\alpha_i))$ and
let
\begin{equation*}
\begin{aligned}
\N(\alpha) & \times^{\reg} \N(l\alpha_i)  = \{(x', x'') \in
\N(\alpha)
\times \N(l \alpha_i) ; \\
& \text{$x'_h$ and $x''_h$ have disjoint spectra for all} \ h\in
\overline{\Omega}_{i}^{\loops} \}.
\end{aligned}
\end{equation*}
By restricting \eqref{E:P1} to $\N(\alpha, l
\alpha_i)_{i,l}:=p_2^{-1}(\N(\alpha + l \alpha_i)_{i,l})$, we obtain
\begin{equation} \label{E:P2}
\N(\alpha)_{i,0} \times^{\reg} \N(l \alpha_i) \overset{p_1}
\longleftarrow \N(\alpha, l \alpha_i)_{i,l} \overset{p_2}
\longrightarrow \N(\alpha + l \alpha_i)_{i,l}.
\end{equation}

\begin{prop}[\cite{KKS08}]\label{P:fiber}
\be[{\rm(a)}]
\item The map $p_2$ in \eqref{E:P2}
is a $GL(\alpha) \times GL(l \alpha_i)$-principal bundle.

\item 
The map $p_1$ in \eqref{E:P2} factors as
$$ \N(\alpha, l\alpha_i)_{i,l} \stackrel{p_1'}{\longrightarrow}
\bigl(\N(\alpha)_{i,0} \times^{\reg} \N(l\alpha_i)\bigr) \times Z(\alpha,
l\alpha_i) \stackrel{p_1''}{\longrightarrow} \N(\alpha)_{i,0}
\times^{\reg} \N(l\alpha_i),$$ where $Z(\alpha, l \alpha_i)$ is the
set of short exact sequences
$$ 0 \longrightarrow V(\alpha) \overset{\phi'} \longrightarrow
V(\alpha + l \alpha_i) \overset{\phi''} \longrightarrow V(l
\alpha_i) \longrightarrow 0,$$  $p_1''$ is the natural projection
and $p_1'$ is an affine fibration.
\ee
\end{prop}

\vspace{.1in}

\begin{cor}[\cite{KKS08}]\label{C:Lag}
\be[{\rm(a)}]
\item For each $\alpha \in Q_{+}$, $\N(\alpha)$ is a Lagrangian
subvariety of $X(\alpha)$.

\item There is a 1-1 correspondence between the set of
irreducible components $\Lambda$ of $\N(\alpha)$
satisfying $\epsor_i(\Lambda)=l$
and those of $\N(\alpha-l\alpha_i)$ satisfying $\epsor_i(\Lambda')=0$.
\ee
\end{cor}

\vspace{.2in}

We denote this 1-1 correspondence by $\Lambda \longmapsto
 \eit^{l}(\Lambda)$.
Set $\mathcal{B} = \coprod_{\alpha \in Q_{+}} \Irr \N(\alpha)$.
Define the maps $\wt: \mathcal{B} \rightarrow -Q_{+} \subset P$,
$\varepsilon_i, \varphi_i : \mathcal{B} \rightarrow \Z \cup
\{-\infty\}$, $\eit, \fit : \mathcal{B} \rightarrow \mathcal{B} \cup
\{0\}$ by
\begin{equation}\label{E:B(infty)}
\begin{aligned}
& \wt(\Lambda)  = -\alpha \ \ \text{for} \ \ \Lambda \in \Irr \N(\alpha), \\
& \varepsilon_i (\Lambda) = \begin{cases} \epsor_i(\Lambda) \
\ & \text{if} \ \ i\in I^{\re}, \\
0 \ \ & \text{if} \ \ i\in I^{\im},
\end{cases}\\
& \varphi_i(\Lambda) = \langle h_i, \wt(\Lambda) \rangle +
\varepsilon_i(\Lambda), \\
& \eit(\Lambda) = \begin{cases} (\eit^{l-1})^{-1} \circ \eit^l
(\Lambda)
 & \text{if}\; \epsor_i(\Lambda)=l>0,\\
 0 & \text{if}\; \epsor_i(\Lambda)=0,
 \end{cases}\\
& \fit(\Lambda)= (\eit^{l+1})^{-1} \circ \eit^l(\Lambda) \qquad
\text{if}\; \epsor_i(\Lambda)=l.
\end{aligned}
\end{equation}

\begin{thm}[\cite{KKS08}] \label{T:B(infty)}
The set $\mathcal{B}$ is a $U_q(\g)$-crystal which is isomorphic to
$B(\infty)$.
\end{thm}


\section{Nakajima's Quiver Variety}

Let $\la \in P^{+}$ be a dominant integral weight and let
$W=W(\lambda)=\bigoplus_{i\in I} W_i$ be an $I$-graded vector space with
$\underline{\wt} (W):= \sum_{i\in I} (\dim W_i) \Lambda_i = \la$.
For each $\alpha \in Q_{+}$, define
$$X(\la;\alpha) = X(\alpha) \oplus \Homi(V(\alpha),W) \oplus \Homi(W,V(\alpha)),$$
where $\Homi(V(\alpha),W) = \bigoplus_{i\in I} \Hom(V_i, W_i)$ and $\Homi(W,V(\alpha))
= \bigoplus_{i\in I} \Hom(W_i, V_i)$. A typical element of
$X(\la;\alpha)$ will be denoted by $(x,t,s)=((x_{h})_{h\in H},
(t_i)_{i\in I}, (s_i)_{i\in I})$. The group $GL(\alpha)$ acts on
$X(\la; \alpha)$ via
$$g \cdot (x,t,s) =((g_{\sink(h)} x_{h} g_{\out(h)}^{-1})_{h\in H},
(t_i g_i^{-1})_{i\in I}, (g_i s_i)_{i\in I}).$$ The symplectic form
$\omega$ on $X(\la; \alpha)$
and the moment map $\mu=(\mu_i: X(\la;\alpha) \rightarrow
\gl(V_i))_{i\in I}$ are given by
\begin{equation*}
\begin{aligned}
\omega((x,t,s), (x',t',s')) & = \sum_{h\in H} \epsilon(h)
\Tr(x_{\oh} x_{h}') + \sum_{i\in I} \Tr(s_i t_i' - s_i' t_i), \\
\mu_i(x,t,s) &= \sum_{\substack{h\in H \\ \out(h)=i}} \epsilon(h)
x_{\oh} x_{h} + s_i t_i.
\end{aligned}
\end{equation*}

For $x=(x_h)_{h\in H} \in X(\alpha)$, an $I$-graded subspace
$U=\bigoplus_{I\in I} U_i$ of $V(\alpha)$ is said to be {\it $x$-stable} if
$x_{h}(U_{\out(h)}) \subset U_{\sink(h)}$ for all $h\in H$.

\begin{defi} \label{D:stable}
A point $(x,t,s) \in X(\la;\alpha)$ is {\it stable} if there is no
nonzero $I$-graded $x$-stable subspace $U=\bigoplus_{i \in I} U_i$
of $V(\alpha)$ such that $t_i(U_i)=0$ for all $i\in I$.
\end{defi}

Let $X(\la;\alpha)^{\st}$ denote the set of all stable points in
$X(\la;\alpha)$. Then the group $GL(\alpha)$ acts freely on $X(\la;
\alpha)^{\st}$ (indeed, if $(x,t,s)$ is stable and $ \Id \neq g \in GL(\alpha)$ satisfies $g \cdot (x,t,s)=(x,t,s)$ then the subspace $\bigoplus_i \text{Im}(g_i-\Id)$ violates the stability condition, see \cite{Nak94}). We define {\it Nakajima's quiver
variety} to be
$$\mathfrak{X}(\la;\alpha) = \mu^{-1}(0) \cap X(\la;\alpha)^{\st}
\big/ GL(\alpha).$$
It is known to be a smooth variety with a symplectic structure induced by $\omega$.
We also set $\N(\lambda;\alpha)= (\N(\alpha) \times \Homi(V(\alpha),W))^{\st}$ and
$$\mathcal{L}(\la;\alpha) = \N(\lambda;\alpha) \big/ GL(\alpha).$$
The definition of the subvariety $\mathcal{L}(\lambda;\alpha)$ is different from the one given in \cite{Nak94} (for quivers without
edge loops), but it yields the same variety (see \cite{Nak94}, Lemma~5.9).

\begin{prop} \label{P:Lag}
For each $\alpha \in Q_{+}$, $\mathcal{L}(\la;\alpha)$ is a
closed Lagrangian subvariety of $\mathfrak{X}(\la;\alpha)$.
\end{prop}
\begin{proof}
By Corollary \ref{C:Lag}, $\N(\alpha)$ is a Lagrangian subvariety of
$X(\alpha)$. Since $\Homi(V,W)$ is clearly a Lagrangian subvariety of
$\Homi(V,W) \oplus \Homi(W,V)$, $\N(\alpha) \times \Homi(V,W)$ is a
Lagrangian subvariety of $X(\la;\alpha)$, which implies $(\N(\alpha)
\times \Homi(V,W))^{\st}$ is a Lagrangian subvariety of
$X(\la;\alpha)^{\st}$. Since $\N(\alpha) \times \Homi(V,W) \subset
\mu^{-1}(0)$, by symplectic reduction, $\mathcal{L}(\la;\alpha)$ is
a Lagrangian subvariety of $\mathfrak{X}(\la;\alpha)$.
\end{proof}

For each $\alpha \in Q_{+}$ and $l\ge 0$, let
\begin{equation*}
\begin{aligned}
E(\la; \alpha, l\alpha_i)  =  \big\{ (& x, x', x'', t,t', \phi',
\phi'') \;; \;  x \in \N(\alpha + l \alpha_i), x' \in \N(\alpha),
x'' \in \N(l
 \alpha_i), \\
& t\in \Homi(V(\alpha+l\alpha_i),W), t' \in \Homi(V(\alpha), W), \\
& 0 \longrightarrow V(\alpha) \overset{\phi'}\longrightarrow
V(\alpha+l\alpha_i) \overset{\phi''} \longrightarrow V(l \alpha_i)
\longrightarrow
0 \ \ \text{is exact}, \\
& \phi' \circ x' = x \circ \phi', \; \phi'' \circ x=x'' \circ
\phi'', t'= t \circ \phi' \big\}
\end{aligned}
\end{equation*}
and consider the canonical projections
\begin{equation}\label{E:P3}
\xymatrix{
E(\la; \alpha, l\alpha_i)\ar[d]^{q_1}\ar[r]^-{q_2}&
\N(\alpha+l\alpha_i) \times \Homi(V(\alpha + l\alpha_i), W)\\
\N(\alpha) \times \Homi(V(\alpha), W) \times \N(l\alpha_i)}
\end{equation}
 given by
\begin{equation*}
\xymatrix{
(x, x', x'', t, t', \phi', \phi'')\ar@{|->}[d]\ar@{|->}[r]&(x, t)\\
(x', t', x'')}
\end{equation*}
It is easy to show that if $(x,t)$ is stable, then $(x',t')$ is also
stable.

Define a function $\epsor_i$ on $\N(\lambda;\alpha)$ by
$\epsor_i(x,t)= \epsor_i(x)$. Note that this function is invariant
under $GL(\alpha)$ and hence descends to
$\mathcal{L}(\lambda;\alpha)$. Set 
\begin{equation*}
\begin{split}
\N(\lambda;\alpha)_{i,l}&=\set{(x,t) \in \N(\lambda;\alpha)}%
{\text{$\epsor_i=l$ on a neighborhood of $(x,t)$}}\\
&=(\N(\alpha)_{i,l} \times \text{Hom}(V(\alpha),W))\cap X(\lambda;\alpha)^{\st},
\end{split}
\end{equation*}
and let
$$\N(\la;\alpha) \times^{\reg} \N(l \alpha_i) =\{(x', t', x'') \in
\N(\la;\alpha) \times \N(l \alpha_i) \,;\, (x',x'') \in \N(\alpha)
\times^{\reg} \N(l \alpha_i) \}.$$ By restricting \eqref{E:P3} to
$\N(\la; \alpha, l \alpha_i)_{i,l}:= q_2^{-1}(\N(\la;\alpha +
l\alpha_i)_{i,l})$, we obtain
\begin{equation} \label{E:P4}
\N(\la; \alpha)_{i,0} \times^{\reg} \N(l \alpha_i) \overset{q_1}
\longleftarrow \N(\la;\alpha, l \alpha_i)_{i,l} \overset{q_2}
\longrightarrow \N(\la; \alpha + l \alpha_i)_{i,l}.
\end{equation}

\begin{prop}\label{P:q-fiber} \hfill
\be[{\rm (i)}]
\item  The map $q_2$ in \eqref{E:P4} is a $GL(\alpha) \times GL(l
\alpha_i)$-principal bundle.

\item
Assume the following conditions:
\eq
&&
\ba{l}
\text{{\rm(a)}\ if $i\in I^\re$, then $\langle h_i,\lambda-\alpha\rangle\ge l$,}\hs{30ex}\\
\text{\rm{(b)}\ if $i\in I^\im$ and $l>0$, then
$\langle h_i,\lambda-\alpha\rangle>0$.}
\ea\label{cond:eps}
\eneq
Then the map $q_1$ in \eqref{E:P4}
 is smooth, locally trivial and with connected
fibers.
\item
If \eqref{cond:eps} is not satisfied, then
$\N(\la; \alpha + l \alpha_i)_{i,l}$ is an empty set.
\ee
\end{prop}
\begin{proof}
(i) For $(x,t) \in \N(\la; \alpha+l\alpha_i)_{i,l}$, we have
$q_2^{-1}(x,t) \simeq p_2^{-1}(x)$ and our assertion follows from
Proposition \ref{P:fiber} (a).

(ii), (iii) Since they are similarly proved as in \cite{Sai} when $i$ is real,
we shall assume that $i$ is imaginary.

We may assume that $l>0$. Assume first that
$\langle h_i,\lambda-\alpha\rangle\le0$.
Set $\alpha = \sum_j k_j \alpha_j \in Q_{+}$.
Then we have
\begin{equation*}
\begin{aligned}
\langle h_i, \la - \alpha \rangle &=
\langle h_i, \la - \sum_j k_j \alpha_j \rangle \\
& = \langle h_i, \la \rangle - k_i a_{ii} -\sum_{j\neq i} k_j a_{ij}
\le 0.
\end{aligned}
\end{equation*}
Then we have $\langle h_i,\lambda\rangle=k_ia_{ii}=k_ja_{ij}=0$.
Hence $W_i=0$ and if there is an arrow $h\cl i\to j$ ($j\not=i$),
then $V(\alpha+l\alpha_i)_j=0$.
Hence $V(\alpha+l\alpha_i)_i\subset \Ker t$ and the
stability condition implies that
$V(\alpha+l\alpha_i)_i=0$, which is a contradiction.

Now we shall assume that $\langle h_i,\lambda-\alpha\rangle>0$.
For $(x',t',x'') \in \N(\la; \alpha)_{i,0} \times^{\reg}
\N(l\alpha_i)$, we have
$$q_1^{-1}(x',t',x'')\simeq p_1^{-1}(x',x'') \times \Homi(V(l\alpha_i),W).$$
Since $q_{1} = p_{1} \times \pi$, where
$$\pi\cl\Homi(V(\alpha+l \alpha_i), W) \times \Homi(V(\alpha), W)
\longrightarrow \Homi(V(\alpha), W)$$ is the natural projection, our
assertion follows from Proposition \ref{P:fiber} (b) once we prove that
for a generic point
$(x,x',x'',t,t',\phi',\phi'')\in E(\lambda;\alpha,l\alpha_i)$
of $q_1^{-1}(x'',t'',x'')$,
then $(x,t,0)$ is a stable point.
Let $U$ be an $I$-grades subspace of $V(\alpha+l\alpha_i)$ such that $t(U)=0$.
Then $U\cap V(\alpha)=0$, and hence $U_j\cap V(\alpha+l\alpha_i)_j=0$
for $j\not=i$. Take $h\in\overline{\Omega}_{i}^{\loops}$.
Then we have a unique $x_h$-invariant decomposition $V(\alpha+l\alpha_i)_i
=V(\alpha)_i\oplus F$. Then $U_i$ is contained in $F$.
We may assume that any eigenvector of $x_h$ in $F$ is not annihilated by $t$.
Hence $t_i(U_i)=0$ implies $U_i=0$.
\end{proof}

As an immediate corollary, there is a 1-1 correspondence
$$\Irr
\N(\la; \alpha)_{i,l} \overset{\sim} \longrightarrow \Irr \N(\la;
\alpha-l\alpha_i)_{i,0},$$
if \eqref{cond:eps} holds. As in \cite[Cor.~3.3]{KKS08}, we deduce by a dimension count that
the irreducible components of $\N(\lambda,\alpha)_{i,l}$ are precisely the intersections of
$\N(\lambda;\alpha)_{i,l}$ with the irreducible components $\Lambda$ of $\N(\lambda;\alpha)$ satisfying
$\epsor_i(\Lambda)=l$. Note also that since $GL(\alpha)$ acts freely on
$\N(\la;\alpha)$, the irreducible components of $\N(\la; \alpha)$ are in 1-1
correspondence with those of $\mathcal{L}(\la; \alpha)$. Hence we
obtain:

\begin{cor}\label{C:1-1} Assume \eqref{cond:eps}. Then
there is a 1-1 correspondence between the irreducible components $\Lambda$
of $\mathcal{L}(\lambda;\alpha+l\alpha_i)$ satisfying $\epsor_i(\Lambda)=l$
and those of $\mathcal{L}(\lambda;\alpha)$ satisfying
$\epsor_i(\Lambda')=0$.
\end{cor}

We denote this 1-1 correspondence by $\Lambda \longmapsto \Lambda'=:
\eit^{l}(\Lambda)$. Observe that
\begin{equation*}
\begin{aligned}
 \Irr \N(\la; \alpha) & = \Irr (\N(\alpha) \times \Homi(V(\alpha),
W))^{\st} \\
&= \{ \tilde{\Lambda} \in \Irr (\N(\alpha) \times \Homi(V(\alpha),
W)) \,;\, \tilde{\Lambda} \cap X(\la; \alpha)^{\st} \neq \emptyset \} \\
&= \{ \Lambda_0 \times \Homi(V(\alpha), W) \,;\, \Lambda_0 \in \Irr
\N(\alpha),\\
& \qquad (\Lambda_0 \times \Homi(V(\alpha), W)) \cap X(\la;
\alpha)^{\st} \neq \emptyset \},
\end{aligned}
\end{equation*}
which defines a map $\psi^{\la}: \Irr \mathcal{L}(\la; \alpha)=\Irr \N(\lambda;\alpha)
\longrightarrow \Irr \N(\alpha)$ given by $\Lambda \mapsto
\Lambda_0$. Note that $\epsor_i(\Lambda) =
\epsor_i(\Lambda_0)$ for all $\Lambda \in \Irr \mathcal{L}(\la;
\alpha)$. Hence, by the definition of $\psi^{\la}$, we obtain the
following commutative diagram:
\begin{equation} \label{E:comm}
\xymatrix{ \Irr \N(\alpha)_{i,l} \ar[r]^-{\eit^{l}} & \Irr \N(\alpha
- l \alpha_i)_{i,0}
\\ \Irr \mathcal{L}(\la; \alpha)_{i,l}
\ar[u]^-{\psi^{\la}} \ar[r]^-{\eit^{l}} & \Irr \mathcal{L}(\la;
\alpha - l \alpha_i)_{i,0}  \ar[u]^{\psi^{\la}} }
\end{equation}

\vskip 1cm


\section{Crystal Structure on $\mathcal{B}^{\la}$}

Let  $\mathcal{B}^{\la} = \coprod_{\alpha \in Q_{+}} \Irr\;
\mathcal{L}(\la; \alpha)$ and for $\Lambda \in \Irr\;
\mathcal{L}(\la; \alpha)$, define
\begin{equation*}
\begin{aligned}
& \wt(\Lambda)  = \la - \alpha, \\
& \varepsilon_i (\Lambda) = \begin{cases} \epsor_i(\Lambda) \
\
& \text{if} \ \ i\in I^{\re}, \\
0 \ \ & \text{if} \ \ i\in I^{\im},
\end{cases}\\
& \varphi_i(\Lambda) = \varepsilon_i(\Lambda) + \langle h_i,
\wt(\Lambda) \rangle.
\end{aligned}
\end{equation*}
If $i\in I^{\re}$, in \cite{Sai}, Saito proved
$$\varphi_i(\Lambda) = \varepsilon_i(\Lambda) + \langle h_i, \la -
\alpha \rangle \ge 0.$$ If $i\in I^{\im}$, write $\alpha = \sum_j k_j
\alpha_j$, then we have $$\varphi_i(\Lambda) = \langle h_i, \la -
\alpha \rangle = \langle h_i, \la \rangle - \sum_j k_j a_{ij} \ge
0.$$ We define the Kashiwara operators $\eit, \fit :
\mathcal{B}^{\la} \rightarrow \mathcal{B}^{\la} \sqcup \{0\}$ by
\begin{equation}\label{E:Eit}
 \begin{aligned}
 \eit(\Lambda) &= \begin{cases} (\eit^{l-1})^{-1} \circ \eit^l (\Lambda)
 & \text{if}\; \epsor_i(\Lambda)=l>0,\\
 0 & \text{if}\; \epsor_i(\Lambda)=0,
 \end{cases}\\
 \fit(\Lambda)&= \begin{cases} (\eit^{l+1})^{-1} \circ \eit^l(\Lambda)
 & \text{if}\; \epsor_i(\Lambda)=l, \, \varphi_i(\Lambda)>0, \\
 0 & \text{if}\; \varphi_i(\Lambda)=0.
 \end{cases}
 \end{aligned}
 \end{equation}
It is straightforward to verify that $\mathcal{B}^{\la}$ is a
$U_q(\g)$-crystal. Moreover, we have

\begin{prop} \label{P:B^la} \hfill

{\rm (a)} The crystal $\mathcal{B}^{\la}$ is connected.

{\rm (b)} If $i\in I^{\im}$ and $\langle h_i, \wt(\Lambda) \rangle
\le -a_{ii}$, then $\eit(\Lambda)=0$.
\end{prop}
\begin{proof}
(a) It suffices to show that if $\epsor_i(\Lambda) =0$ for all
$i\in I$, then $\alpha =0$ and $\Lambda = \{0\}$, which was already
proved in \cite{KKS08}.

(b) If $\Lambda'\seteq\eit(\Lambda)\not=0$, 
then we have $\langle h_i,\wt(\Lambda')\rangle
=\langle h_i,\wt(\Lambda)\rangle+a_{ii}\le0$ and hence
$\fit(\Lambda')=0$ by Proposition \ref{P:q-fiber} (iii).
Hence it is a contradiction.
%
%
\end{proof}

\vskip 3mm

Define a map $\Psi^{\la}: \mathcal{B}^{\la} \longrightarrow
\mathcal{B} \ot T_{\la} \ot C$ by $\Lambda \longmapsto
\psi^{\la}(\Lambda) \ot t_{\la} \ot c$.

\begin{thm} \label{T:main}
The map $\Psi^{\la}$ is a strict crystal embedding.
\end{thm}
\begin{proof}
If $\Lambda \in \Irr \mathcal{L}(\la; \alpha)$ with $\alpha \in
Q_{+}$, then $\psi^{\la}(\Lambda)=\Lambda_0 \in \Irr \N(\alpha)$ and
we have
$$\wt(\Psi^{\la}(\La)) = \wt(\Lambda_0 \ot t_{\la} \ot c) = \la - \alpha = \wt(\Lambda).$$

If $i \in I^{\re}$, by the definition of tensor product of crystals,
we have
\begin{equation*}
\begin{aligned}
\varepsilon_i(\Lambda_0 \ot t_{\la} \ot c) &=
\max(\varepsilon_i(\Lambda_0), - \langle h_i, \la - \alpha \rangle
), \\
\varphi_i(\Lambda_0 \ot t_{\la} \ot c) &= \max(\varphi_i(\Lambda_0)
+ \langle h_i, \la \rangle, 0).
\end{aligned}
\end{equation*}
Since $\varepsilon_i(\Lambda) =\varepsilon_i( \Lambda_0)$, we have
\begin{equation*}
\varphi_i(\Lambda_0) = \langle h_i, \la \rangle =
\varepsilon_i(\Lambda_0) + \langle h_i, \la - \alpha \rangle   =
\varepsilon_i(\Lambda) + \langle h_i, \la - \alpha \rangle =
\varphi_i (\Lambda) \ge 0.
\end{equation*}
Hence we obtain
$$\varepsilon_i(\Lambda_0 \ot t_{\la} \ot c) = \varepsilon_i(\Lambda),
\quad \varphi_i(\Lambda_0 \ot t_{\la} \ot c) = \varphi_i(\Lambda).$$
If $i \in I^{\im}$, then $$\varepsilon_i(\Lambda) = 0 =
\varepsilon_i(\Lambda_0 \ot t_{\la} \ot c), \quad
\varphi_i(\Lambda) = \langle h_i, \la - \alpha \rangle =
\varphi_i(\Lambda_0 \ot t_{\la} \ot c).$$

It remains to show that $\Psi^{\la}$ commutes with $\eit, \fit$
$(i\in I)$. By Example \ref{exam:C}, we have
\begin{equation*}
\begin{aligned}
\fit(\Lambda_0 \ot t_{\la} \ot c) &= \begin{cases} \fit \La_0 \ot
t_{\la} \ot c & \text{if} \ \varphi_i(\La) >0, \\
0 & \text{if} \ \varphi_i(\La)=0,
\end{cases} \\
\eit(\Lambda_0 \ot t_{\la} \ot c) & = \begin{cases} \eit(\La_0) \ot
t_{\la} \ot c & \text{if} \ i \in I^{\re}, \ \varphi_i(\La) \ge 0, \\
\eit(\La_0) \ot t_{\la} \ot c & \text{if} \ i \in I^{\im}, \ \langle
h_i, \la - \alpha \rangle + a_{ii} >0, \\
0 & \text{if} \ i \in I^{\im}, \ \langle h_i, \la - \alpha \rangle +
a_{ii} \le 0.
\end{cases}
\end{aligned}
\end{equation*}

If $\varphi_i(\La)=0$, then $\fit(\La)=0$ and hence $\Psi^{\la}(\fit
\La)=0 = \fit(\La_0 \ot t_{\la} \ot c)$. If $\varphi_i(\La)>0$, then
using the commutative diagram \ref{E:comm}, we obtain
\begin{equation*}
\begin{aligned}
\fit \Psi^{\la} (\La) & = \fit (\La_0 \ot t_{\la} \ot c) = \fit
\La_0 \ot t_{\la} \ot c \\
&= (\eit^{l+1})^{-1} \circ \eit^{l} (\La_0) \ot t_{\la} \ot c =
\psi^{\la} ((\eit^{l+1})^{-1} \circ \eit^{l} (\La)) \ot t_{\la} \ot
c \\
& = \psi^{\la}(\fit \La) \ot t_{\la} \ot c = \Psi^{\la}(\fit \La).
\end{aligned}
\end{equation*}

Note that $\eit \La =0$ if and only if $\eit \La_0 =0$. Hence if
$i\in I^{\re}$ and $\varphi_i(\La) \ge 0$, using the commutative
diagram \ref{E:comm}, we have
\begin{equation*}
\begin{aligned}
\eit \Psi^{\la} (\La) & = \eit (\La_0 \ot t_{\la} \ot c) = \eit
\La_0 \ot t_{\la} \ot c \\
&= (\eit^{l-1})^{-1} \circ \eit^{l} (\La_0) \ot t_{\la} \ot c =
\psi^{\la} ((\eit^{l-1})^{-1} \circ \eit^{l} (\La)) \ot t_{\la} \ot
c \\
& = \psi^{\la}(\eit \La) \ot t_{\la} \ot c = \Psi^{\la}(\eit \La).
\end{aligned}
\end{equation*}
Similarly, if $i\in I^{\im}$ and $\langle h_i, \la - \alpha \rangle +
a_{ii} >0$, one can verify $\eit \Psi^{\la} (\La) = \Psi^{\la} (\eit
\La)$. Finally, if $i \in I^{\im}$ and $\langle h_i, \la - \alpha
\rangle + a_{ii} \le 0$, by Proposition \ref{P:B^la}, we have $\eit
(\La)=0$ and hence $\Psi^{\la}(\eit \La)=0 = \eit (\La_0 \ot t_{\la}
\ot c)$, which completes the proof.
\end{proof}

As a corollary we obtain the geometric realization of the crystal
$B(\la)$.

\begin{cor} \label{C:main}
The crystal $\mathcal{B}^{\la}$ is isomorphic to the highest weight
crystal $B(\la)$.
\end{cor}
\begin{proof}
Let $\mathbf{1}_{\la}$ be the unique element of $\mathcal{B}^{\la}$
satisfying $\varepsilon_i(\mathbf{1}_{\la})=0$ for all $i\in I$.
Then $\mathbf{1}:= \psi^{\la} (\mathbf{1}_{\la})$ is the unique
element of $\mathcal{B}$ such that $\varepsilon_i(\mathbf{1})=0$ for
all $i\in I$ and we have $\Psi^{\la}(\mathbf{1}_{\la}) = \mathbf{1}
\ot t_{\la} \ot c$. Hence $\mathcal{B}^{\la}$ is isomorphic to the
connected component of $\mathcal{B} \ot T_{\la} \ot C$ containing
$\mathbf{1} \ot t_{\la} \ot c$. Since $\mathcal{B} \cong B(\infty)$,
by Proposition \ref{P:B(la)}, we conclude $\mathcal{B}^{\la} \cong
B(\la)$.
\end{proof}

\vskip 10mm

\vskip 1cm

\providecommand{\bysame}{\leavevmode\hbox
to3em{\hrulefill}\thinspace}

\end{document}
